\documentclass[12pt]{article}

\usepackage{amsmath}
\usepackage{amsfonts}
\usepackage{graphics}
\usepackage{bbm,authblk}
\usepackage{epstopdf}
\usepackage{graphicx}

\usepackage[final]{listings}
\usepackage[numbers,sort&compress]{natbib}

\newtheorem{theorem}{Theorem}[section]

\begin{document}

\title{Numerical computation of the Rosenblatt distribution and applications}

\author{{Nikolai N. Leonenko\textsuperscript{a}\thanks{CONTACT N.~N. Leonenko. Email: leonenkon@cardiff.ac.uk} ~and Andrey Pepelyshev\textsuperscript{a}} }
\affil{ \textsuperscript{a}School of Mathematics, Cardiff University, Cardiff, UK}







\maketitle

\begin{abstract}
The Rosenblatt distribution plays a key role in the limit theorems for non-linear functionals
of stationary Gaussian processes with long-range dependence.
We derive new expressions for the characteristic function of the Rosenblatt distribution.
Also we present a novel accurate approximation of all eigenvalues of the Riesz integral operator
associated with the correlation function of the Gaussian process and propose an efficient algorithm for computation of the density of the Rosenblatt distribution.
We perform Monte-Carlo simulation for small sample sizes to demonstrate the appearance of the Rosenblatt distribution for several functionals of stationary Gaussian processes with long-range dependence.
\end{abstract}




\section{Introduction}

The phenomenon of long-range dependence (also called long memory) is one of exciting area of research in the probability theory for last few decades due to
non-standard normalizations and non-Gaussian limiting distributions of nonlinear functionals \cite{pipiras2017long}.
The Rosenblatt distribution serves a significant role in the study of  this phenomenon
which occurs in economics, finance, hydrology, turbulence, cosmology and physics \cite{doukhan2002theory,beran2017statistics}.

The Rosenblatt distribution was introduced in \cite{rosenblatt1961independence}
and later it was investigated by many researchers, see \cite{taqqu1975weak,taqqu1979convergence,dobrushin1979non} among many others.
The Rosenblatt distribution appears as the limiting distribution of some popular functionals,
see \cite{doukhan2002theory,berman1979high,rosenblatt1979some} and references therein.

The known analytical form of the characteristic function of the Rosenblatt distribution contains a series which converges in a neighbourhood of zero, see
\cite{rosenblatt1961independence, albin1998note,veillette2013properties}.
It turns out that the Rosenblatt distribution is infinitely divisible \cite{veillette2013properties} and self-decomposable, moreover, it belongs to the Thorin class \cite{maejima2013distribution,leonenko2017non,leonenko2017rosenblatt}.

The density of the Rosenblatt distribution exists and bounded, however, the closed analytical form of the density is unknown, see \cite{veillette2013properties,maejima2013distribution,leonenko2017non,leonenko2017rosenblatt}.
The Edgeworth expansion was used to approximate the density of the Rosenblatt distribution \cite{veillette2012berry}. The numerical evaluation of the density of the Rosenblatt distribution was proposed in \cite{veillette2013properties}.

The present paper is organized as follows.
In Section 2 we review the known facts on the Rosenblatt distribution.
In Section 3 we provide the novel analytic forms of the characteristic function of the Rosenblatt distribution on the entire line and a compelling viewpoint on the structure of the Rosenblatt distribution.
In Section 4 we propose a accurate approximation of all eigenvalues of the Riesz integral operator, which allows us to perform fast numerical computation of the density of the Rosenblatt distribution.

In Section 5
we propose a time-efficient algorithm for simulation of the stationary Gaussian sequences with long-range dependence with the power correlation function and the Mittag-Lefler correlation function.
Furthermore, we demonstrate small sample properties of four popular functionals for analysis of sequences of correlated random variables including their non-standard normalization and convergence to the Rosenblatt distribution. Specifically, we consider  estimation of the mean in Section \ref{sec:5m},  the correlation function in Section \ref{sec:5c},
the sojourn functionals in  Section \ref{sec:5s} and path roughness in Section \ref{sec:5r}.
All these problems are essentials of the statistical analysis of stationary Gaussian sequences with  long-range dependence.

\section{Formal statement}

We consider the Rosenblatt distribution with shape parameter $a$, zero mean,
unit variance and the characteristic function
\begin{equation}\label{eq:Rchfun}
 \phi(z)=\exp\left( \frac{1}{2} \sum_{k=2}
 (2iz)^k\frac{c_{a,k} \sigma_a^k }{k} \right),~a\in[0,1/2],~z\in S_0\subset\mathbb{R},
\end{equation}
where $S_0$ is a small neighbourhood of zero,
$$\sigma_a=\sqrt{(1-2a)(1-a)/2}$$
and
$$
 c_{a,k}=\int_0^1\cdots\int_0^1  |x_1-x_2|^{-a} |x_2-x_3|^{-a}
 \cdots |x_{k-1}-x_k|^{-a} |x_k-x_1|^{-a} dx_1\cdots dx_k.
$$
In the rest of this section we describe the Rosenblatt distribution following  \cite{veillette2013properties}.
The random variable $V$ from the Rosenblatt distribution can be given as
\begin{equation}\label{eq:V}
V=\sum_{n=1}^\infty \lambda_{a,n} (\varepsilon_n^2-1),
\end{equation}
where $\varepsilon_n$ are i.i.d. random variables from the standard normal distribution
and
$\lambda_{a,1},\lambda_{a,2},\ldots$  are eigenvalues of the Riesz integral operator
$\tilde K_a: L^2(0,1)\to L^2(0,1)$ defined as
$$
 (\tilde K_af)(x)=\sigma_a \int_0^1 |x-u|^{-a}f(u)du.
$$
The representation \eqref{eq:V} means that $V$ is a specific instance of second-order Wiener chaos \cite{nourdin2012convergence}.
These eigenvalues  satisfy the relation
$$
\sum_{n=1}^\infty \lambda_{a,n}^k=\begin{cases}
\infty,&k=1,\\
c_{a,k} \sigma_a^k,&k=2,3,\ldots.
\end{cases}
$$
Moreover, we have
 \begin{eqnarray*}
\sum_{n=1}^\infty \lambda_{a,n}^2&=&1/2, \\
\sum_{n=1}^\infty \lambda_{a,n}^3&=&\frac{2\sigma_a^3}{(1-a)(2-3a)}\beta(1-a,1-a),
 \end{eqnarray*}
where $\beta(u,v)=\int_0^1x^{u-1}(1-x)^{v-1}dx$ is the beta function.
The Rosenblatt distribution is infinitely divisible with
the L\'evy  density
$$
m(x)=\frac{1}{2x}\sum_{n=1}^\infty
\exp\left(-\frac{x}{2\lambda_{a,n}}\right),~x>0,
$$
that shows self-decomposability of the Rosenblatt distribution \cite{veillette2013properties}.

\section{Main results}

The Laplace transform of $V$ defined in \eqref{eq:V} is given by
$$
 \phi_{LT}(s)=\mathbb{E}(e^{-sV})=\exp\left( -\sum_ {n=1}^\infty
 \left(\frac{1}{2}\ln(1+2\lambda_{a,n}s )-\lambda_{a,n}s\right)  \right),
 ~s>-\frac{1}{2\lambda_{a,1}},
$$
see \cite{albin1998note}. We note that the formal expansion of the sum in the above expression gives
$$
 -\frac{1}{2}\sum_ {n=1}^\infty
 \ln(1+2\lambda_{a,n}s ) +
 s\sum_ {n=1}^\infty\lambda_{a,n},
$$
where the latter summand equals infinity.
Using the classical Taylor expansion of the logarithm
$$
\ln(x)=\sum_{k = 1}^{\infty} \frac{(-1)^{k-1}}{k} (x-1)^k,~ x\in(0,2),
$$  which does not converge for all positive $x$, we obtain
\begin{equation}
 \phi_{LT}(s)=\exp\left( \frac{1}{2}\sum_{k=2}^\infty \frac{(-2s)^k}{k}
 \left( \sum_ {n=1}^\infty
 \lambda_{a,n}^k\right)  \right),~|s|<\frac{1}{2\lambda_{a,1}},
 \label{eq:LTstandard}
\end{equation}
which is related to \eqref{eq:Rchfun} as $\phi(z)=\phi_{LT}(-iz)$, $z\in\mathbb{R}$.
We note that the expressions of the characteristic function \eqref{eq:Rchfun}  and the Laplace transform \eqref{eq:LTstandard} cannot be used for numerical calculation because they are defined at a neighbourhood of zero.

In next theorem, we derive new expressions using
 three other expansions of  the logarithm which converge on the full domain $(0,\infty)$ and one integral representation. In other words, we construct analytic continuations for the expressions \eqref{eq:Rchfun} and~\eqref{eq:LTstandard}.

\begin{theorem} \label{th1}
 The characteristic function $\phi(z)$ of the Rosenblatt distribution is given by $\phi(z)=\phi_{LT}(-iz)$,
  where the Laplace transform $\phi_{LT}(s)$ with $s>-\frac{1}{2\lambda_{a,1}}$ admits
  \begin{itemize}
  \item[(i)] the representation using the domain-scaled expansion
 \begin{eqnarray}\nonumber
 \ln ( \phi_{LT}(s))
 = \sum_ {n=1}^\infty \left(
 \frac{(\lambda_{a,n}s)^2}{1+\lambda_{a,n}s}
  - \sum_{k = 2}^{\infty} \frac{1}{2k-1}\left(\frac{\lambda_{a,n}s}{1+\lambda_{a,n}s}\right)^{2k-1}
 \right),\label{eq:LTds}
\end{eqnarray}
 \item[(ii)] the representation using the Ramanujan expansion
\begin{eqnarray}
 \ln ( \phi_{LT}(s))= \sum_ {n=1}^\infty \lambda_{a,n}s  \ln(1+2\lambda_{a,n}s)
 \sum_{k=1}^\infty \frac{1}{2^{k}(1+(1+2\lambda_{a,n}s)^{2^{-k}})},
\end{eqnarray}
  \item[(iii)] the representation using the Ramanujan-Bradley expansion
\begin{eqnarray}
 \ln ( \phi_{LT}(s))=\frac{1}{2} \sum_ {n=1}^\infty \sum_{k=1}^\infty 2^{k-1}
 \left((1+2\lambda_{a,n}s)^{2^{-k}}-1\right)^2,
\end{eqnarray}
  \item[(iv)] the representation using the integral form
  $$
     \ln ( \phi_{LT}(s))=\int_{0}^{ s} \sum_ {n=1}^\infty \frac{2\lambda_{a,n}^2u}{2\lambda_{a,n} u+1}du.
  $$
  \end{itemize}
\end{theorem}

The proof of Theorem \ref{th1} is deferred to the online supplement.
We note that  the evaluation of the characteristic function
of the Rosenblatt distribution requires computation of the infinite sums with a special care because
numerical computation of $\sum_{i=1}^\infty  \lambda_{a,n}^2$ is problematic for $a\in(0.35,0.5)$ due to very slow convergence of the series $\sum_{i=1}^\infty n^{2a-2}$.
To overcome difficulties of numerical computation, we propose to
find smallest integer $M_\epsilon$ such that
$$
 \sum_{n=M_\epsilon+1}^\infty  \lambda_{a,n}^k<\epsilon
$$
 for $k=3,4,\ldots$, where
 $\epsilon$ is a small positive number, for example, $\epsilon=0.0001$.
Then we can write the random variable $V$ defined in \eqref{eq:V} in the form $V=W+Z$, where
$$
 W=\sum_{n=1}^{M_\epsilon}  \lambda_{a,n} (\varepsilon_n^2-1),~
 Z=\sum_{n=M_\epsilon+1}^\infty  \lambda_{a,n} (\varepsilon_n^2-1).
$$
Due to the choice of $M_\epsilon$, the random variable $Z$ has approximately
the normal distribution with mean zero and variance
$$\sigma^2_\epsilon=1-2\sum_{n=1}^{M_\epsilon}  \lambda_{a,n}^2
=2\sum_{n=M_\epsilon+1}^\infty  \lambda_{a,n}^2.$$
Indeed, the characteristic function of $Z$ has the form
\begin{eqnarray*}
 \phi_Z(z)&=&\exp\left( -\sum_ {n=M_\epsilon+1}^\infty
 \left(\frac{1}{2}\ln(1-2\lambda_{a,n}iz )+\lambda_{a,n}iz\right)  \right)\\
 &\stackrel{|z|<\frac{1}{2\lambda_{a,M_\epsilon}}}{=}&\exp\left( \frac{1}{2}\sum_{k=2}^\infty \frac{(2iz)^k}{k}
 \left( \sum_ {n=M_\epsilon+1}^\infty
 \lambda_{a,n}^k\right)  \right)\\
 &\approx&\exp\left(-z^2  \left( \sum_ {n=M_\epsilon+1}^\infty
 \lambda_{a,n}^2\right) +\frac{4}{3} i^3z^3 \epsilon \right)\\
 &\approx& \exp\left(-z^2 \sigma^2_\epsilon/2 \right).
\end{eqnarray*}
Moreover, we have
$$
\max_{|z|\le\frac{1}{2\lambda_{a,M_\epsilon}}}
\Big|\phi_Z(z)-\exp\left(-z^2 \sigma^2_\epsilon/2 \right)\Big|\to 0 \text{~as~}\epsilon\to0.
$$

This argument provides a clear constructive view on the Rosenblatt distribution.
Specifically, the random variable from the Rosenblatt distribution can be simulated as
$$
  V_\epsilon= \sigma_\epsilon\varepsilon_0+ \sum_{n=1}^{M_\epsilon}  \lambda_{a,n} (\varepsilon_n^2-1),
$$
where $\varepsilon_n$ are i.i.d. random variables from the standard normal distribution.
The density of the Rosenblatt distribution can be computed by the inverse Fourier transform
of the characteristic function
$$
 \phi_\epsilon(z)=\exp\left( -z^2 \sigma^2_\epsilon/2-\sum_ {n=1}^{M_\epsilon}
 \left(\frac{1}{2}\ln(1-2\lambda_{a,n}iz )+\lambda_{a,n}iz\right)  \right).
$$
Applying the Stein method for the random variable $V_\epsilon$ and taking the limit as $\epsilon\to0$, we obtain the following characterizing identity,
which follows from \cite[Ch. 3]{arras2019stein}.
\begin{theorem}
 Let $X$ be a random variable with zero mean. The Stein equation
 \begin{equation}
  \mathbb{E}Xf(X)=\mathbb{E}\int_0^\infty \Big(f(X+x)-f(X)\Big)
  \frac{1}{2}\sum_{n=1}^\infty
\exp\left(-\frac{x}{2\lambda_{a,n}}\right)dx
\label{eq:stein}
 \end{equation}
holds for all bounded Lipschitz function $f(\cdot)$ if and only if $X$ has the Rosenblatt distribution with parameter~$a$.
\end{theorem}
Let us use the Stein equation \eqref{eq:stein} with $f(x)=x$,
we refer to \cite[Ch. 3]{arras2019stein} for a discussion on a class of suitable functions.
Then the right hand side  of the Stein equation is
$$
 \mathbb{E}\int_0^\infty (X+x-X)\frac{1}{2}\sum_{n=1}^\infty
\exp\left(-\frac{x}{2\lambda_{a,n}}\right)dx
=
 2\sum_{n=1}^\infty  \lambda_{a,n}^2=1
$$
and
the left hand side of the Stein equation is the second moment, which equals one from the definition of the Rosenblatt distribution. Thus, the Stein equation  for $f(x)=x$ and $X=V$ becomes $\mathbb{E}V^2=1$.
Taking $f(x)=x^2$ and $X=V$ in  the Stein equation \eqref{eq:stein} we obtain
$$
 \mathbb{E}V^3=
  \mathbb{E}\int_0^\infty (2Vx+x^2)\frac{1}{2}\sum_{n=1}^\infty
\exp\left(-\frac{x}{2\lambda_{a,n}}\right)dx
=
 8\sum_{n=1}^\infty  \lambda_{a,n}^3.
$$
For $f(x)=x^3$ and $X=V$ in  the Stein equation \eqref{eq:stein} we obtain
    \begin{equation*}
 \mathbb{E}V^4
 =
  \mathbb{E}\int_0^\infty (3V^2x-3Vx+x^3)\frac{1}{2}\sum_{n=1}^\infty
\exp\left(-\frac{x}{2\lambda_{a,n}}\right)dx
=
 48\sum_{n=1}^\infty  \lambda_{a,n}^4+3.
    \end{equation*}
Using the first four moments (mean, variance, skewness, and kurtosis), the Rosenblatt distribution can be  approximated by the Pearson family of distributions.

\section{Computational aspects}

Although moments of the Rosenblatt distribution depend on the multidimensional integrals $c_{a,k}$,
$$
 \mathbb{E}\,V^k=\frac{1}{i^k}\left.\frac{d^k}{dz^k}\phi(z)\right|_{z=0}=
 \begin{cases}
 0,&k=1,\\
 1,&k=2,\\
 8\sigma_a^3 c_{a,3},&k=3\\
 48 \sigma_a^4 c_{a,4}+12 \sigma_a^4 c_{a,2}^2,&k=4\\
 \vdots&k=5,6,\ldots,
 \end{cases}
$$
simulation and computation of the density of the Rosenblatt distribution
is not possible without the eigenvalues $\lambda_{a,n}$.
It was shown in \cite{dostanic1998spectral} that
the eigenvalues $\lambda_{a,n}$  have the asymptotic behaviour
$\lambda_{a,n}=C_a n^{a-1}(1+o(n^{-\delta}))$ as $n\to\infty$ for any $\delta\in(0,1)$ and
admit the accurate approximation
$$
 \lambda_{a,n}\cong C_a n^{a-1},~ n>20,
$$
 where
$$
C_a=\frac{2\sigma_a}{\pi^{1-a}}\Gamma(1-a)\sin(\pi a/2),
$$
see \cite{veillette2013properties} for details.
From extensive numerical calculation of eigenvalues for various $a\in(0,0.5)$, we obtain that
the eigenvalues $\lambda_{a,n}$ admit the accurate approximation
\begin{equation} \label{eq:app-eigs}
 \lambda_{a,n}\cong \begin{cases}
 (1+0.1409a)\sqrt{\pi^a\Gamma(1-a)}\sqrt{\frac{1}{2}-a},&n=1,\\
 C_a n^{a-1}+\frac{5}{4} a^{1.05}\sqrt{\Gamma(a+\frac{1}{2})-1}\,n^{a-2.2},&n=2,3,\ldots.
 \end{cases}
\end{equation}

\begin{figure}[h]
\begin{center}
\includegraphics[width=0.39\linewidth]{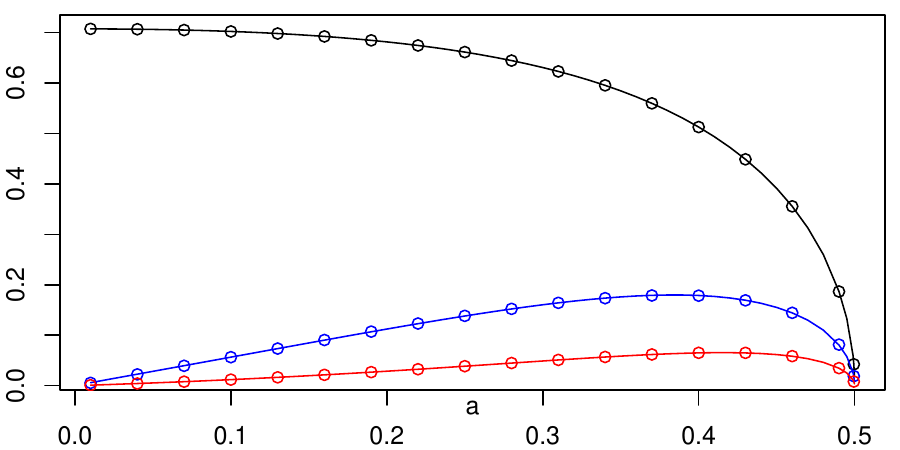}
\includegraphics[width=0.39\linewidth]{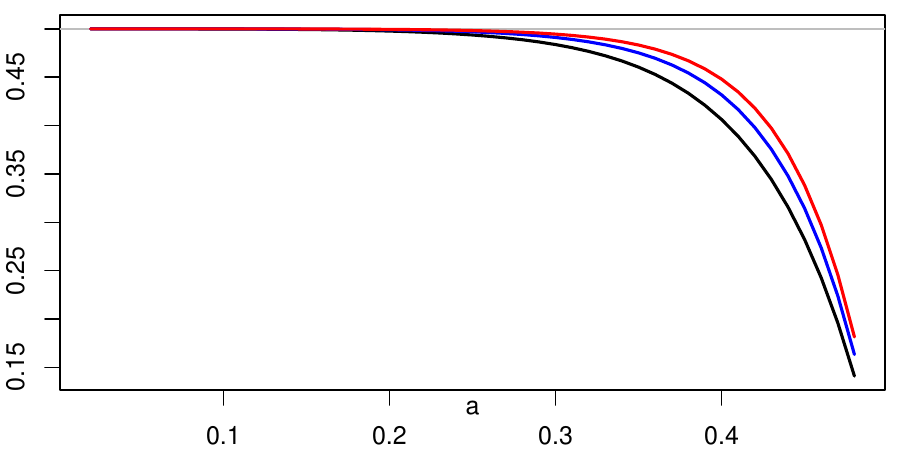}
\end{center}
\caption{Left: Eigenvalues $\lambda_{a,1}$ (black), $\lambda_{a,2}$ (blue) and $\lambda_{a,8}$ (red) which are computed numerically (circles) and via approximation \eqref{eq:app-eigs} (solid line) for $a\in(0,0.5)$.
Right: The sum $\sum_{n=1}^M \lambda_{a,n}^2$ as a function of $a$ for $M=100$ (black),
$M=500$ (blue) and $M=2000$ (red).}
\label{fig:eigs}
\end{figure}

In Figure \ref{fig:eigs} we demonstrate that eigenvalues $\lambda_{a,n}$ for various $a$ and $n$ computed via approximation \eqref{eq:app-eigs} are close to values computed numerically
using the algorithm described in \cite{veillette2013properties}.
Figure~\ref{fig:eigs} also shows the behaviour
of $S_{a,M}=\sum_{n=1}^M \lambda_{a,n}^2$
which show the contribution of $W$. We note that $\sigma^2_\epsilon=1-2S_{a,M}$ describes
the contribution of the normal component $Z$.
We note that for  $a\approx0.45$ we have
the situation, where both components $W$ and $Z$ have approximately equal contribution
for making the shape of the Rosenblatt distribution, that is,
$S_{a,M}\approx\sigma^2_\epsilon$ for  $a\approx0.45$.

Let us find $M_\epsilon$ from the condition
$\sum_{n=M_\epsilon+1}^\infty  \lambda_{a,n}^3\cong\epsilon$.
We can see from Table \ref{tab:Me} that the case of $a\cong0.44$ requires the largest value of $M_\epsilon$. This means that the accurate computation of the Rosenblatt distribution requires a big number of eigenvalues for $a\approx0.44$.

\begin{table}[!hhh]
\caption{Values of $M_\epsilon$ for various $\epsilon$ and $a$.}
\begin{center}
\begin{tabular}{|l|r|r|r|r|r|r|r|}
  \hline
  &$a=0.1$ &$a=0.2$ &$a=0.3$&$a=0.35$ &$a=0.4$&$a=0.44$&$a=0.48$\\
  \hline
  $M_\epsilon|_{\epsilon=10^{-3}}$\rule{0mm}{6mm}&2&3& 7& 13& 24&34 & 13\\
  $M_\epsilon|_{\epsilon=10^{-4}}$\rule{0mm}{6mm}&3&9&48&133&409&909&630\\
\hline
\end{tabular}
\end{center}
\label{tab:Me}
\end{table}

\begin{figure}[h]
\begin{center}
\includegraphics[width=0.39\linewidth]{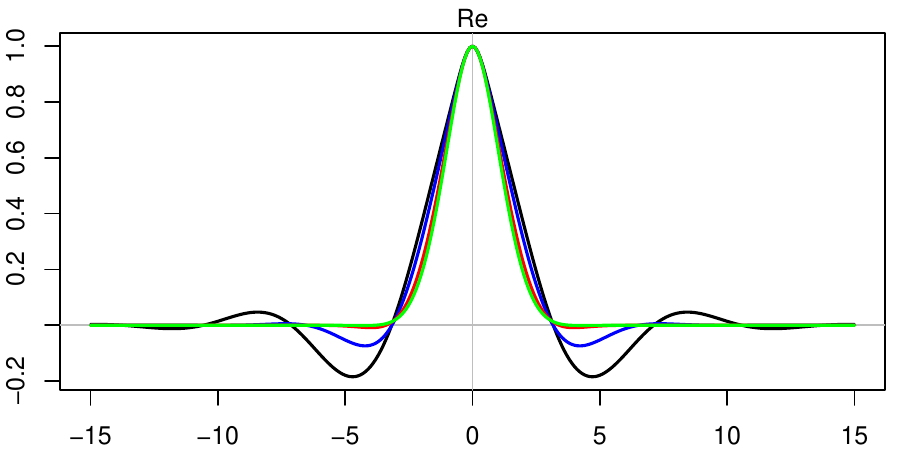}
\includegraphics[width=0.39\linewidth]{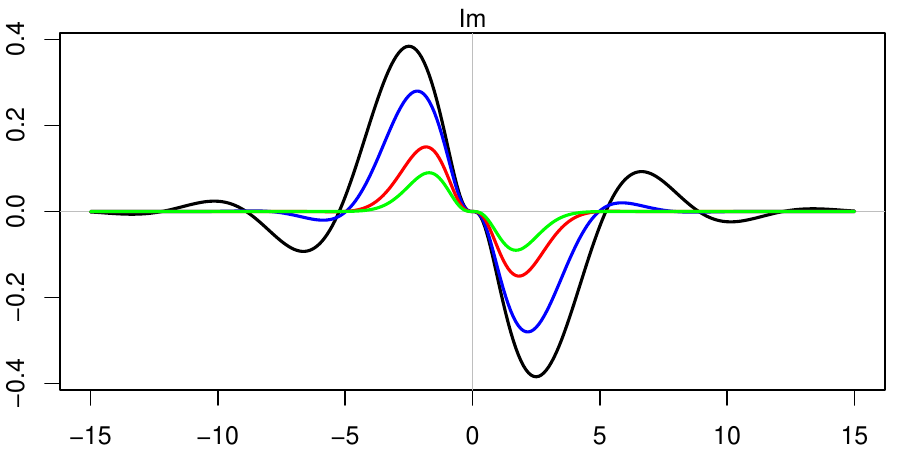}
\end{center}
\caption{Real part (left) and imaginary part of the characteristic function $\phi_\epsilon(z)$
with $M_\epsilon=100$ for $a=0.2$ (black), $a=0.3$ (blue), $a=0.4$ (red)
and $a=0.44$ (green).}
\label{fig:chf}
\end{figure}

In Figure \ref{fig:chf} we depict the characteristic function $\phi_\epsilon(z)$ for various $a$.
We can see that the shape of $\phi_\epsilon(z)$ tends to the characteristic function of the normal distribution as $a\to0.5.$

\begin{figure}[h]
\begin{center}
\includegraphics[width=0.33\linewidth]{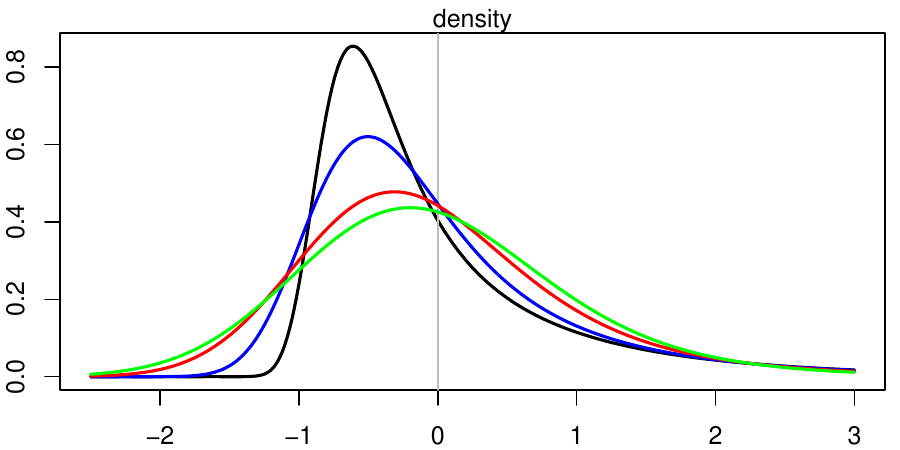}
\includegraphics[width=0.33\linewidth]{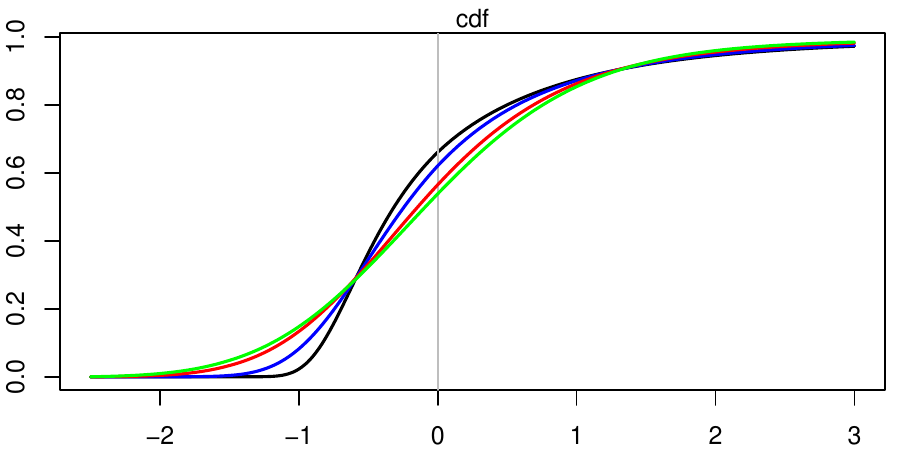}
\includegraphics[width=0.33\linewidth]{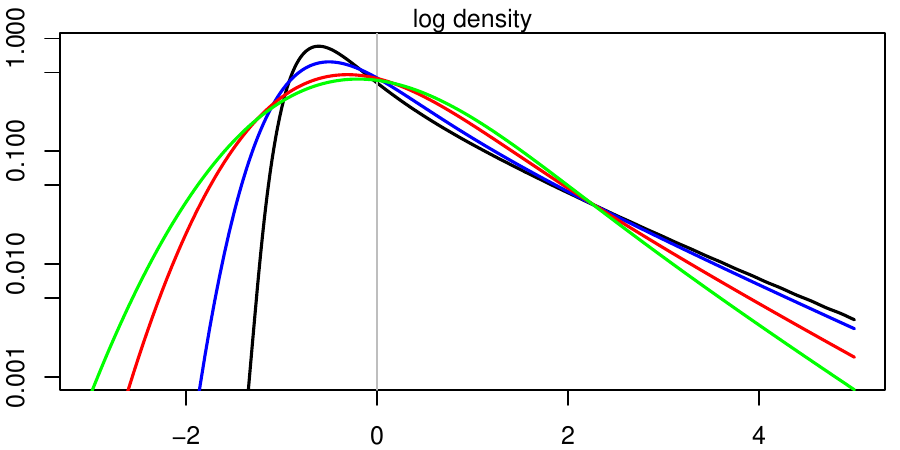}
\end{center}
\caption{The density (left), the logarithm of the density (bottom) and the cdf (right) of the Rosenblatt distribution which is computed via the inverse Fourier transform of the characteristic function $\phi_\epsilon(z)$
with $M_\epsilon=100$ for $a=0.2$ (black), $a=0.3$ (blue), $a=0.4$ (red)
and $a=0.44$ (green).}
\label{fig:dens}
\end{figure}

In Figure \ref{fig:dens} we depict the density, the logarithm of the density and the cumulative distribution function (cdf) of the Rosenblatt distribution for various $a$. We can see that the density of Rosenblatt distribution is close to the density of shifted chi-square distribution
for $a\in(0,0.2)$ and close to the density of the normal distribution for $a\in(0.4,0.5)$.

\section{Numerical study with applications of the Rosenblatt distribution}

The Rosenblatt distribution with parameter $a$ appears as a limiting distribution of several functionals for stationary Gaussian sequences with the correlation function of the form
$$
 r(t)=\frac{l(|t|)}{|t|^a},
$$
where $l(\cdot)$ is a slowly varying function, that implies long-range dependence if $a\in(0,1)$.
Typical examples of such correlation function  are
$r(t)=(1+t^2)^{-a/2}$, $r(t)=(1+|t|)^{-a}$ and $r(t)=E_a(-|t|^a)$, where $E_a(\cdot)$ is the Mittag-Leffler function.
For running our numerical study,
we firstly propose an efficient algorithm for simulation of  long sequences with LRD,
see \cite{bardet2003generators} for review of simulation algorithms, which were found to be very time consuming in our settings.

\subsection{Simulation of long Gaussian sequences}

The traditional way of simulation of a stationary Gaussian sequence $X_1,\ldots,X_n$ with zero mean, unit variance and the specified correlation function is via simulation of the vector from the multivariate normal distribution with the use of the Cholesky decomposition or the eigenvalue decomposition of the covariance matrix.
However, this way is computationally infeasible if the length of Gaussian sequence is larger than 10000 due to the problem with the decomposition of the covariance matrix of large size.

The second way of simulation is to approximate the Gaussian sequence by the autoregressive process. This approach works well when the specified correlation function is close to zero for large lags and is not suitable for simulation of sequences with LRD.

The third way of simulation is based on the approximation
$$
 X_j\approx\sum_{k=1}^M \sqrt{b_k} X_j^{(k)},~j=1,\ldots,n,
$$
where $X_1^{(k)},\ldots,X_n^{(k)}$ is a Gaussian autoregressive sequence of order 1
with zero mean, unit variance and the correlation function $e^{-\lambda_k |t|}$
and the specified correlation function $r(t)$ allows the approximation in the form
\begin{equation}
 r(t)\approx\sum_{k=1}^M b_k e^{-\lambda_k |t|}.
 \label{eq:appr-r(t)}
\end{equation}
We note that the coefficients $b_k$ could be computed via the minimization of an error of the approximation \eqref{eq:appr-r(t)} with some $\lambda_k$, for example, $\lambda_k=e^{-(k-1)}$. To avoid this high-dimensional optimization, we propose the following methodology for finding $b_k$ and $\lambda_k$ in the approximation~\eqref{eq:appr-r(t)}.

Suppose that the correlation function $r(t)$ has the representation
$$
 r(t)=\int_0^\infty e^{-tx}p(x)dx,
$$
where $p(x)$ is a density.
Then we can construct the approximation \eqref{eq:appr-r(t)} with $b_k=1/M$ and
$\lambda_k$ be random values from the distribution with the density $p(x)$.
Alternatively, $\lambda_k$ can be chosen as the $(k-0.5)/M$-quantile of the distribution with the density $p(x)$.

Let us take a decreasing positive sequence $q_0,q_1,q_2,\ldots$ such that
$1=q_0>q_1>q_2>\ldots>0$ and $\lim_{k\to\infty} q_k=0$.
Define $\tau_k$ as the $q_k$-quantile of the distribution with density $p(x)$, that is,
$$
 q_k=\int_0^{\tau_k} p(x)dx,~k=0,1,2,\ldots.
$$
We note that  $\tau_1,\tau_2,\ldots$ is a decreasing sequence,  $\tau_1>\tau_2>\ldots$ and $\lim_{k\to\infty} \tau_k=0$.
Then by splitting the integration domain $(0,\infty)$ with breakpoints $\tau_1,\tau_2,\ldots$ we
obtain that
$$
 r(t)=
 \sum_{k=1}^\infty \int_{\tau_{k}}^{\tau_{k-1}} e^{-tx}p(x)dx
$$
and, consequently, the correlation function $r(t)$ has the approximation \eqref{eq:appr-r(t)} with
$$
 b_k=q_{k-1}-q_{k}\mbox{~and~some~}
 \lambda_k\in(\tau_{k-1},\tau_{k}),
$$
that follows from
$$
 \int_{\tau_{k}}^{\tau_{k-1}} e^{-tx}p(x)dx
 \approx e^{-t \lambda_k} \int_{\tau_{k}}^{\tau_{k-1}} p(x)dx
 =(q_{k-1}-q_{k}) e^{-t \lambda_k}.
$$
We recommend to take $\lambda_1=\tau_1$ and
$$
 \lambda_k=\sqrt{\tau_{k-1}\tau_k},~k=2,3,\ldots.
$$
The approximation \eqref{eq:appr-r(t)} is accurate if $M$ is large and $\max \{b_1,\ldots,b_M\}$ is small.

\subsection{Simulation of long Gaussian sequences with LRD}

Suppose that
\begin{equation}\label{RV_measure}
   p(x)=\frac 1 {\Gamma(\alpha)}x^{\alpha-1}l(1/x), \;\alpha>0, ~x>0,~ x\to 0,
\end{equation}
where $l(\cdot)$ is a slowly varying function,
which describes the behaviour of the density $p(x)$ at zero.
It follows from the Tauberian-Abelian theorem that the property \eqref{RV_measure} is equivalent to
\begin{equation}\label{RV_corr}
    r(t)=\frac {l(|t|)}{|t|^{\alpha}}, \; t\to\infty,
\end{equation}
which describes the behaviour of the correlation function $r(t)$ at infinity, see e.g. \cite{feller1958introduction}. We note that LRD occurs if $\alpha\in(0,1).$

For simulation of sequences  with $r(t)$ in the form \eqref{RV_corr},
it was proposed in \cite{leonenko2005convergence} to take
$$b_k=c/k^{(1+a)}, ~\lambda_k=1/k,$$
where $c$ is a constant such that $\sum_{k=1}^M b_k=1$.
This choice of $\lambda_k$ is not convenient in numerical studies due to slow convergence to zero.

Let us develop an approximation for the correlation function $r(t)=1/(1+|t|)^a$.
From \cite{BNL} we have that
$$
 \frac{1}{(1+|t|)^{a}}=\int_0^\infty e^{-tx}p(x)dx,~a\in(0,1),
$$
where
$$
 p(x)=\frac{1}{\Gamma(a)}x^{a-1}e^{-x}.
$$
We recommend to choose the  sequence $q_k$ to be fast-decreasing.
Specifically, we propose to take
\begin{equation}
 q_k=\begin{cases}
 0.98,&k=1,\\
 0.9\gamma^{k-2},&k=2,3,\ldots,
 \end{cases}
 ~\gamma=e^{-(2-a)a},
 ~M=\left\lceil\frac{2}{a}\right\rceil+8,
 \label{eq:bk-lamk}
\end{equation}
The above choice of $q_k$ provides the reasonable accuracy of the approximation \eqref{eq:appr-r(t)} even with small~$M$. The simulation algorithm is deferred to the online supplement.

\begin{figure}[h]
\begin{center}
\includegraphics[width=0.56\linewidth]{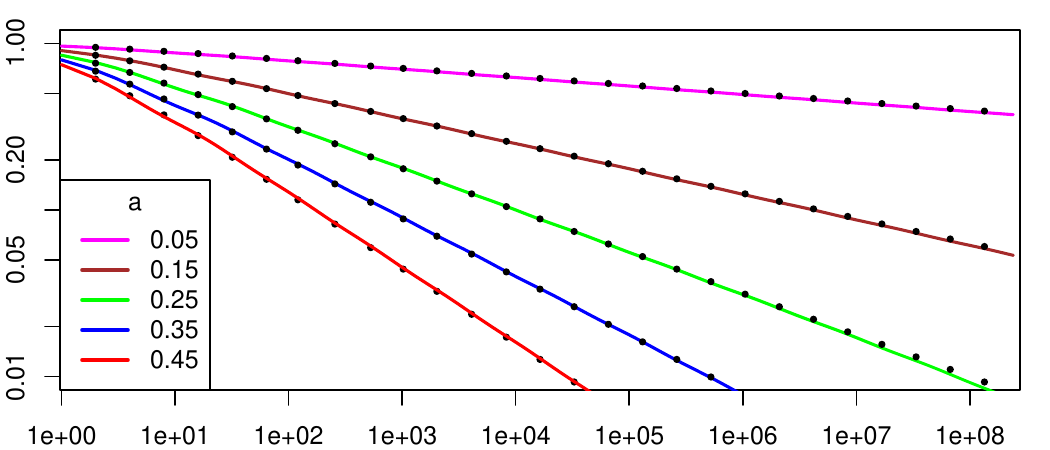}
\end{center}
\caption{The approximation \eqref{eq:appr-r(t)} with parameters \eqref{eq:bk-lamk} (solid line) of the correlation function $r(t)=1/(1+|t|)^a$ (dotted line) for various $a$ in log-log scale. }
\label{fig:cor-appr}
\end{figure}

In Figure \ref{fig:cor-appr} we demonstrate the good accuracy of the approximation \eqref{eq:appr-r(t)} with parameters \eqref{eq:bk-lamk}  of the correlation function $r(t)=1/(1+|t|)^a$  for various $a\in(0,0.5).$

\begin{figure}[h]
\begin{center}
\includegraphics[width=0.46\linewidth]{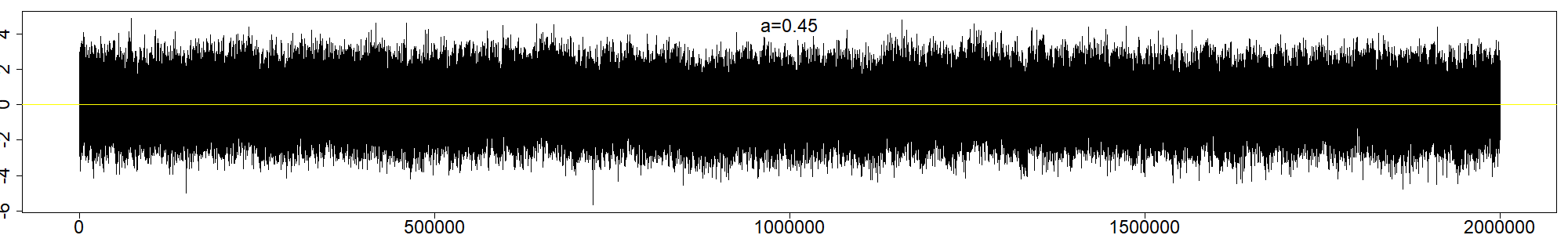}
\includegraphics[width=0.46\linewidth]{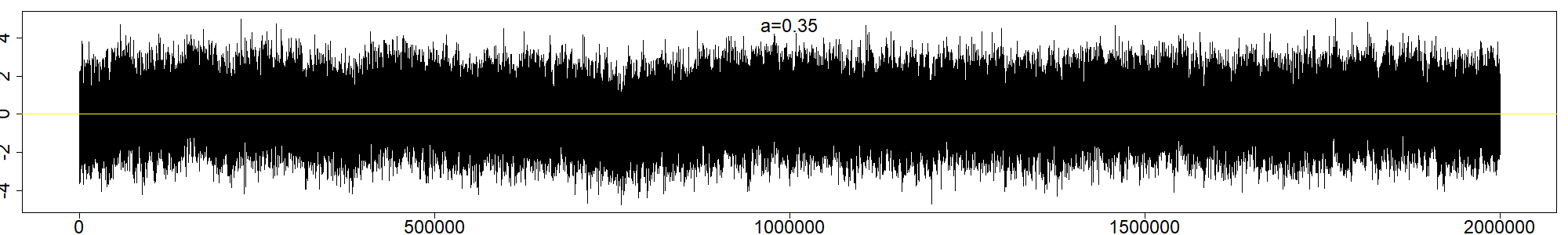}
\includegraphics[width=0.46\linewidth]{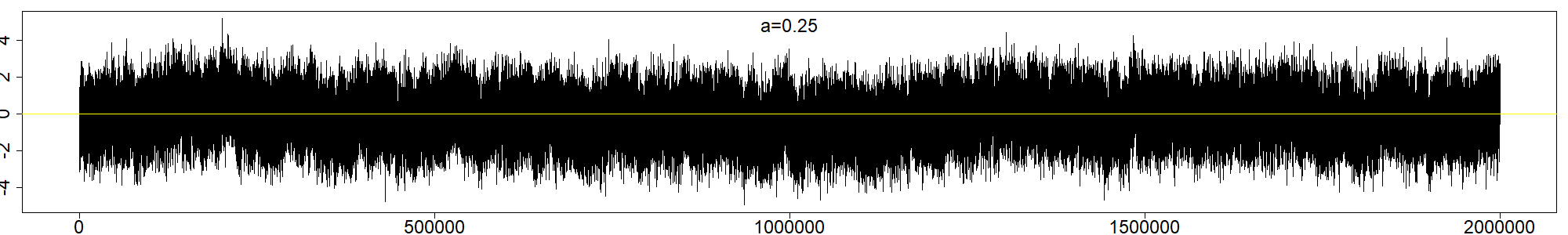}
\includegraphics[width=0.46\linewidth]{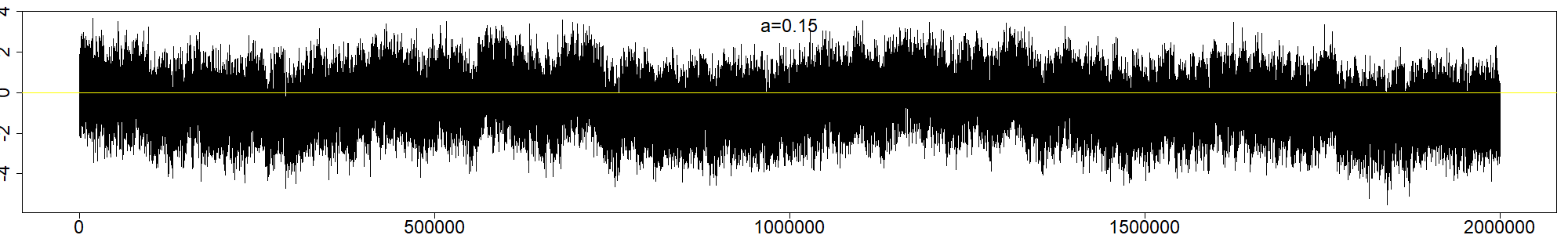}
\end{center}
\caption{The realizations of a stationary Gaussian sequence with zero mean, unit variance and the correlation function $r(t)=1/(1+|t|)^a$ for $a=0.15$, $0.25$, $0.35$ and $0.45$. }
\label{fig:seq}
\end{figure}

In Figure \ref{fig:seq} we depict several realizations of a stationary Gaussian sequence with zero mean, unit variance and the correlation function $r(t)=1/(1+|t|)^a$ to illustrate the phenomenon of long-range dependence.
We can see that the deviation of a local trend of realizations from zero increases as $a$ decreases.
We note that a sequence of length 2000000 is needed for rather good accuracy of estimation of the parameter $a$ if $a\approx0.25$.

\begin{figure}[h]
\begin{center}
\includegraphics[width=0.56\linewidth]{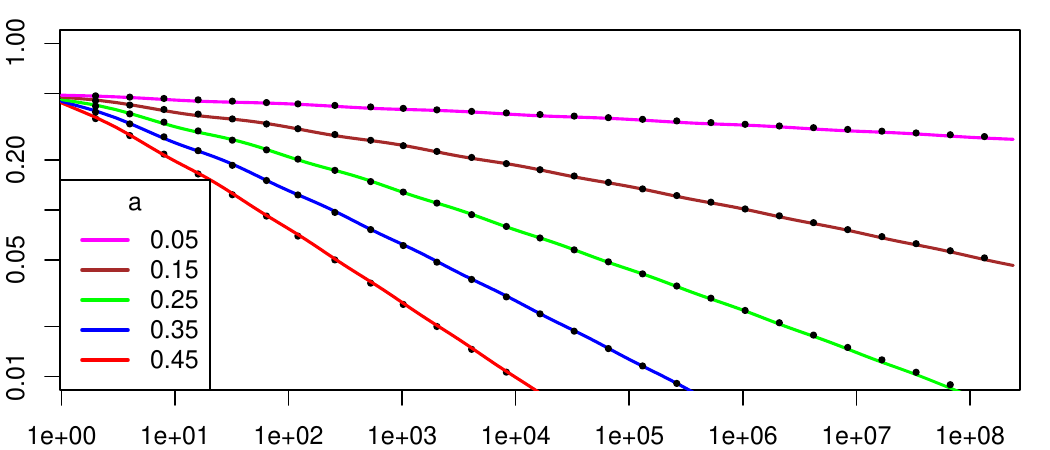}
\end{center}
\caption{The approximation \eqref{eq:appr-r(t)} with parameters \eqref{eq:bk-lamk-ml} (solid line) of the Mittag-Leffler correlation function $r(t)=E_a(-|t|^a)$ (dotted line) for various $a$ in log-log scale. }
\label{fig:cor-appr-ml}
\end{figure}

As  a second example, we develop an  approximation for the correlation function
$r(t)=E_a(-|t|^a)$, where $E_a(t)$ is the Mittag-Lefler function with parameter $a\in(0,1)$.
From \cite{BNL} we have that
$$
 E_a(-|t|^a)=\int_0^\infty e^{-tx}p_{ml}(x)dx,~a\in(0,1),
$$
where
$$
 p_{ml}(x)=\frac{\sin(a\pi)}{\pi}\frac{x^{a-1}}{1+2\cos(a\pi) x^a+x^{2a}},~x>0,
$$
is the density of the Lamperti distribution,
which has the quantile function
$$
 Q_{ml}(u)=\left(\frac{\sin(u a\pi)}{\sin((1-u) a \pi)}\right)^{1/a},~u\in(0,1),
$$
see \cite{james2010lamperti}.
We recommend to choose the  sequence $q_k$ in the form
\begin{equation}
 q_k=\begin{cases}
 0.98,&k=1,\\
 0.9,&k=2,\\
 0.7,&k=3,\\
 0.5\gamma^{k-4},&k=4,5,\ldots,
 \end{cases}
 ~\gamma=e^{-(2-a)a},
 ~M=\left\lceil\frac{2}{a}\right\rceil+8,
 \label{eq:bk-lamk-ml}
\end{equation}
The above choice of $q_k$ for the Mittag-Lefler correlation function provides the reasonable accuracy of the approximation \eqref{eq:appr-r(t)} even with small~$M$.
In Figure \ref{fig:cor-appr-ml} we demonstrate the good accuracy of the approximation \eqref{eq:appr-r(t)} with parameters \eqref{eq:bk-lamk-ml}  of the correlation function
$r(t)=E_a(-|t|^a)$  for various $a\in(0,0.5).$

\subsection{Estimation of the mean for a sequence of special structure}
\label{sec:5m}
Consider a stationary sequence $Y_1,Y_2,\ldots$, which is given by
$$
 Y_k=\theta+ \sigma H_2(E_k),
$$
where $H_2(x)=x^2-1$ and $E_1,E_2,\ldots$ is a stationary Gaussian sequence
with zero mean, unit variance and correlation function $r(k)$. The sequence $\{E_k\}$ can be viewed as discretization of a stationary Gaussian process and the sequence $\{H_2(E_k)\}$ is a stationary sequence with zero mean and variance $2$. Consider the estimator of $\theta$ in the form
$$
 \hat \theta_n=\frac{1}{n}\sum_{k=1}^n Y_k.
$$
The estimator $\hat\theta_n$ has asymptotically
the scaled Rosenblatt distribution if the sequence $\{E_k\}$ has the correlation function $r(t)=\frac{l(|t|)}{|t|^a}$ with $a\in(0,0.5)$ as $t\to\infty$
and has  the normal distribution
otherwise, see \cite{ivanov2002asymptotic,leonenko2006weak}.

Let us make a numerical study to obtain the distribution of the estimator for small sample sizes.
Specifically, we consider the random variable
$$
 Z_n=\frac{\sigma_a}{ n^{1-a}}\sum_{k=1}^n H_2(E_k),~a\in(0,0.5),
$$
which
has asymptotically the Rosenblatt distribution with parameter $a\in(0,0.5).$

\begin{figure}[h]
\begin{center}
\includegraphics[width=0.2\linewidth]{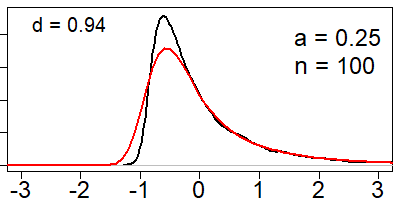}
\includegraphics[width=0.2\linewidth]{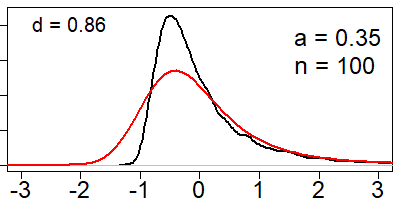}
\includegraphics[width=0.2\linewidth]{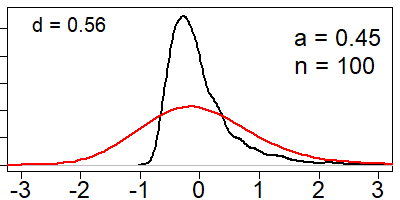}\\
\includegraphics[width=0.2\linewidth]{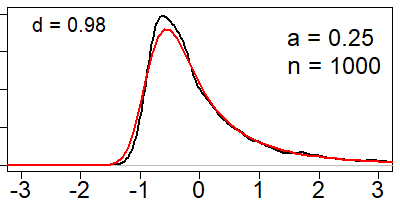}
\includegraphics[width=0.2\linewidth]{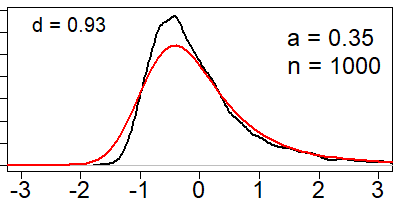}
\includegraphics[width=0.2\linewidth]{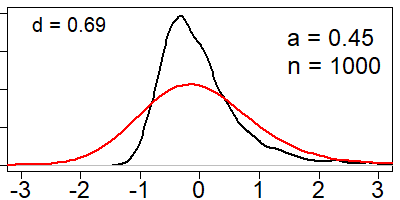}\\
\includegraphics[width=0.2\linewidth]{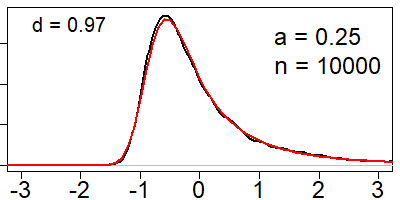}
\includegraphics[width=0.2\linewidth]{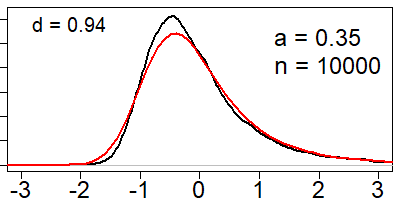}
\includegraphics[width=0.2\linewidth]{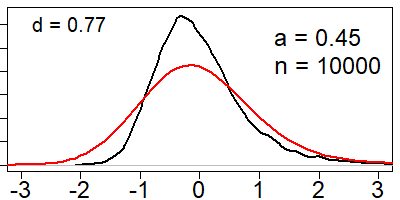}\\
\includegraphics[width=0.2\linewidth]{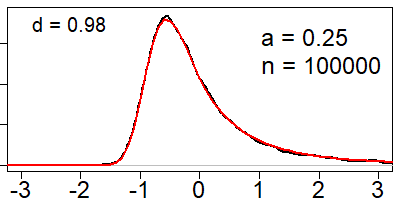}
\includegraphics[width=0.2\linewidth]{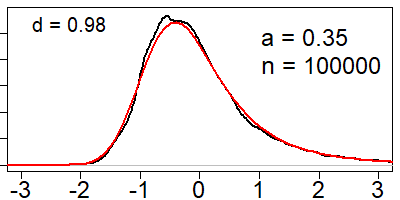}
\includegraphics[width=0.2\linewidth]{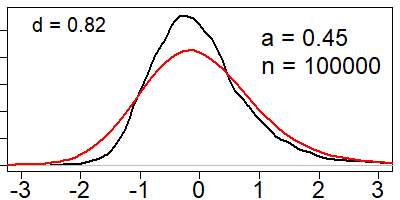}
\end{center}
\caption{The empirical density of $Z_n$ (black) and the density of the Rosenblatt distribution (red) for various $n$ and $a=0.25$, $0.35$ and $0.45$;
$d$ is the empirical standard deviation of $Z_n$. }
\label{fig:md}
\end{figure}

In Figure \ref{fig:md} we can see that the empirical density of $Z_n$
for $a=0.25$ and various $n$ is very close to
the density of the Rosenblatt distribution and this holds for all $a\in(0,0.25)$.
However, such proximity is not observed for larger values of $a$.
Nevertheless, the empirical density of $Z_n$ for $a=0.45$ becomes more close to
the density of the Rosenblatt distribution as $n$ increases.
We would like to note that  the empirical density of $Z_n$ for $a=0.45$ and small $n$
is close to the Rosenblatt distribution with some small value of $a$.

\subsection{Estimation of the correlation function for a stationary Gaussian sequence with LRD}
\label{sec:5c}

Let $E_1,E_2,\ldots$ be a stationary Gaussian sequence with zero mean, unit variance and
correlation function $r(k)=\mathrm{Cov}(E_s,E_{s+k})=\mathbb{E}(E_sE_{s+k})$.
Suppose that the correlation function $r(k)$ has the shape
$$
 r(k)=\frac{l(k)}{k^a}.
$$
As well known, the classical estimator of $r(k)$ is given by
$$
\hat r(k)=\frac{1}{n}\sum_{j=1}^{n-k} E_jE_{j+k},
$$
which has asymptotically
the scaled Rosenblatt distribution if $a\in(0,0.5)$ and
the normal distribution otherwise, see \cite{rosenblatt1961independence,rosenblatt1979some}.
The fact that the asymptotic distribution of $\hat r(k)$ does not depend on $k$ for $a\in(0,0.5)$ is very remarkable.
Specifically, the limiting distribution of
$$
  R_{k,n}=\frac{\sigma_a}{n^{1-a}}\sum_{j=1}^{n-k}\Big(E_jE_{j+k}-r(k)\Big),
  ~a\in(0,0.5),~k=0,1,\ldots.
$$
is given by the Rosenblatt distribution.

Let us make a numerical study to obtain the distribution of $R_{k,n}$ for small sample sizes.
We note that $R_{0,n}=Z_n$ and, therefore, the empirical distribution of $R_{0,n}$
is shown in Figure \ref{fig:md}.

\begin{figure}[h]
\begin{center}
\includegraphics[width=0.2\linewidth]{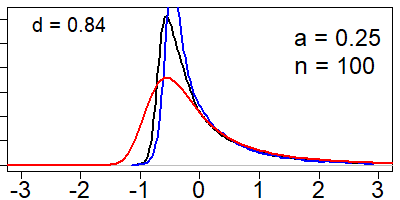}
\includegraphics[width=0.2\linewidth]{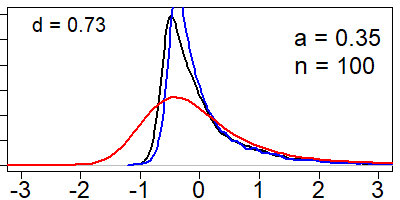}
\includegraphics[width=0.2\linewidth]{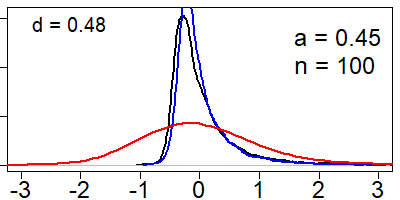}\\
\includegraphics[width=0.2\linewidth]{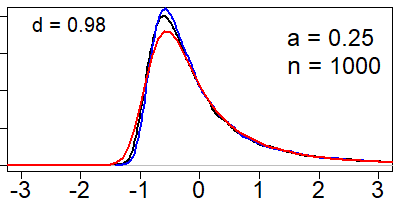}
\includegraphics[width=0.2\linewidth]{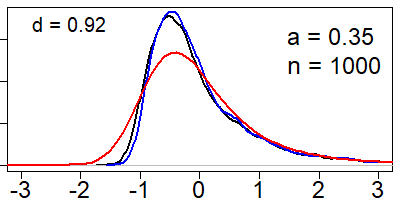}
\includegraphics[width=0.2\linewidth]{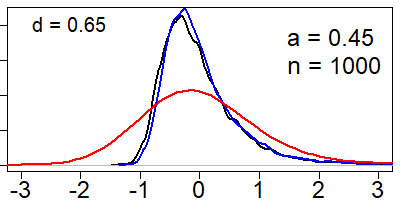}\\
\includegraphics[width=0.2\linewidth]{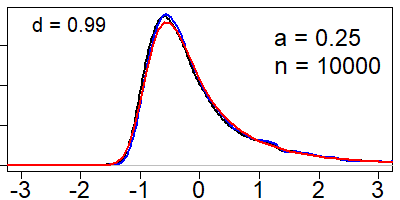}
\includegraphics[width=0.2\linewidth]{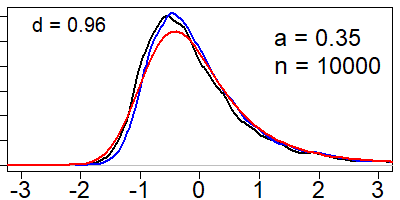}
\includegraphics[width=0.2\linewidth]{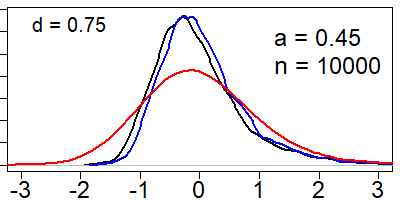}
\end{center}
\caption{The empirical density of $R_{10,n}$ (black) and $R_{20,n}$ (blue) and the density of the Rosenblatt distribution (red) for various $n$ and $a=0.25$, $0.35$ and $0.45$;
$d$ is the empirical standard deviation of $R_{10,n}$. }
\label{fig:cor}
\end{figure}

In Figure \ref{fig:cor} we can see that the empirical distribution of $R_{k,n}$ for small $n$ almost does not depends on $k$, that confirms the theoretical statement.
In particular, the empirical distribution of $R_{k,n}$ is similar to the empirical distribution of $Z_{n}$ for various $a$ and $n$.

\subsection{Estimation of sojourn functionals}
\label{sec:5s}

Let $E(t)$ be a stationary Gaussian process with zero mean, unit variance and
 correlation function $r(t)=\frac{l(|t|)}{|t|^a}$. Consider the sojourn functional
$$
 M_t(u)=\int_0^t \mathbf{1}_{\{|E(s)|>u\}} (s) ds,
$$
where $u>0$ and
$$
  \mathbf{1}_{\{|E(s)|>u\}}(s)=\begin{cases}
  1,&|E(s)|>u,\\
  0,&|E(s)|\le u,
  \end{cases}
$$
is the indicator function.
We interpret $M_t(u)$ as the time spent by the process $|E(s)|$ above the level $u$ for $s\in[0,t]$.
Following \cite{berman1979high}, the expansion of $M_t(u)$ is the form
$$
 M_t(u)=2(1-\Phi(u))+u\phi(u)\int_0^t(E^2(s)-1)ds+2\phi(u)\sum_{j=2}^\infty
 \frac{H_{2j-1}}{2j!}\int_0^t H_{2j}(E(s))ds,
$$
where $H_j(\cdot)$ is the $j$-th Hermit polynomial, and
$$
 \mathrm{Var}(M_t(u))=4\int_0^t (t-s)\int_0^{r(s)}(\phi(u,u;q)-\phi(u,u;-q))dqds,
$$
where
$$\phi(u,v;q)=\frac{1}{2\pi\sqrt{1-q^2}}\exp\left(  -\frac{x^2-2quv+v^2}{2(1-q^2)} \right)$$
is the standard bivariate normal density with correlation $q$.
From \cite{berman1979high} we also have that
the functional
$$
 \frac{M_t(u)-2t(1-\Phi(u))}{2u\phi(u)\sqrt{\int_0^t (t-s)r^2(s)ds}}
$$
has asymptotically the Rosenblatt distribution.

Let us make a numerical study to obtain the distribution of
$$
 S_{u,n}=\frac{\sigma_a}{ n^{1-a} }\cdot\frac{1}{u\phi(u)}\left(
 \sum_{j=1}^n \mathbf{1}_{\{|E_j|>u\}} (j)- 2n(1-\Phi(u))
 \right)
$$ for small sample sizes, where $E_1,\ldots,E_n$ is a stationary Gaussian sequence with zero mean and
correlation function $r(t)=1/(1+|t|)^a$.

\begin{figure}[h]
\begin{center}
\includegraphics[width=0.2\linewidth]{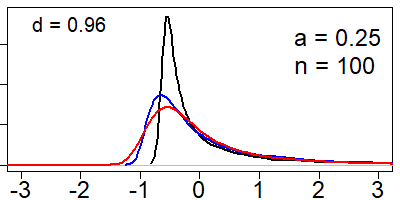}
\includegraphics[width=0.2\linewidth]{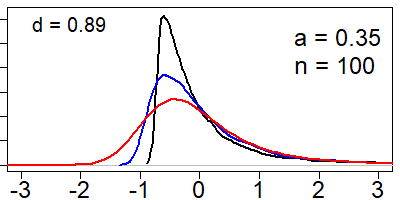}
\includegraphics[width=0.2\linewidth]{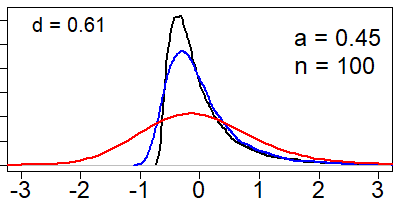}\\
\includegraphics[width=0.2\linewidth]{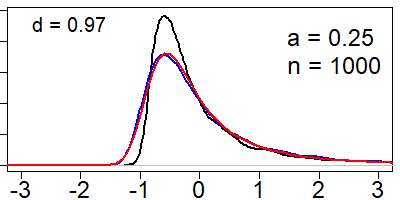}
\includegraphics[width=0.2\linewidth]{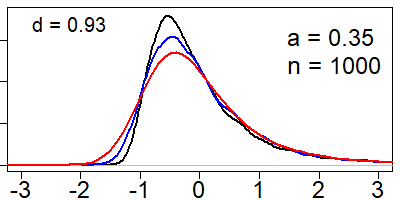}
\includegraphics[width=0.2\linewidth]{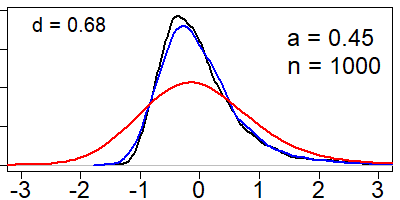}\\
\includegraphics[width=0.2\linewidth]{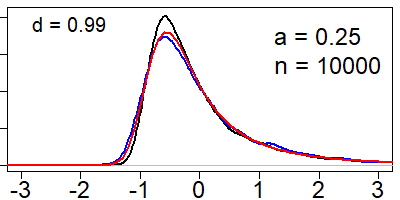}
\includegraphics[width=0.2\linewidth]{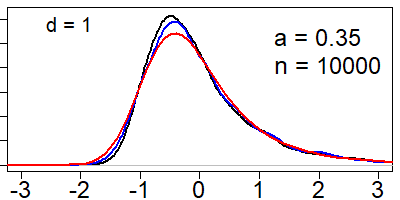}
\includegraphics[width=0.2\linewidth]{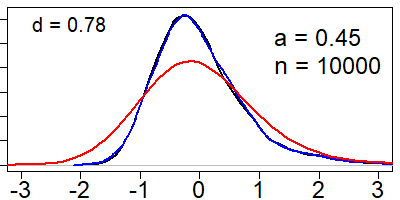}
\end{center}
\caption{The empirical density of $S_{2,n}$ (black) and $S_{1.5,n}$ (blue) and the density of the Rosenblatt distribution (red) for various $n$ and $a=0.25$, $0.35$ and $0.45$;
$d$ is the empirical standard deviation of $S_{2,n}$. }
\label{fig:soj}
\end{figure}

In Figure \ref{fig:soj} we can see that the empirical distribution of $S_{u,n}$ for small $n$ depends slightly on $u$.
In particular, the empirical distribution of $S_{u,n}$ tends to the Rosenblatt distribution as
$n$ increases.

\subsection{Roughness of the fBm path}
\label{sec:5r}

Let $X(t)$ be a stochastic process, $t\in[0,1]$.
Define the quadratic variation of $X(t)$  by
$$
 V_n=\sum_{j=1}^n \Big(X(\textstyle{\frac{j}{n}}) -X(\textstyle{\frac{j-1}{n}})\Big)^2,
$$
which can be interpreted as the path roughness and often used in finance and geophysics.

Suppose that $X(t)$ is the fractional Brownian motion (fBm) with Hurst parameter $H\in(0,1)$
and covariance function $r(t,s)=(s^{2H}+t^{2H}-|t-s|^{2H})/2$.
Then the distribution of
$$
F_n=\frac{V_n-\mathbb{E}(V_n)}{\sqrt{\mathrm{Var}(V_n)}}.
$$
is asymptotically normal when $H\in(0,3/4]$
and converges to the Rosenblatt distribution with parameter $a=2-2H$ when
$H\in(3/4,1)$ as $n$ increases, see
\cite{dobrushin1979non,nourdin2012selected,nourdin2012convergence}.
The quadratic variation $V_n$ is useful in practice because the classical estimator of $H$ is based on the fact that $n^{2H-1}V_n\to1$ in probability as $n\to\infty$.

Let us make a numerical study to obtain the distribution of
$$
 G_n=\frac{\sigma_a}{ n^{1-a} }\cdot \frac{n^{2-a}}{1.04-1.5a}\sum_{j=1}^n \left[
 \Big(X(\textstyle{\frac{j}{n}}) -X(\textstyle{\frac{j-1}{n}})\Big)^2-n^{a-2}\right]
$$
for small sample sizes, where $X(t)$ is the fractional Brownian motion with $H=1-a/2$,
which can be simulated using the FFT-based algorithm from \cite{wood1994simulation}
implemented in the R package {\it SuperGauss}.

\begin{figure}[h]
\begin{center}
\includegraphics[width=0.2\linewidth]{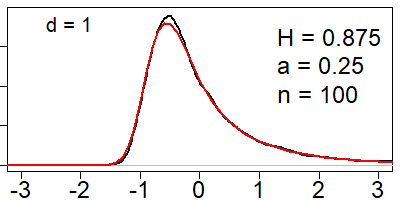}
\includegraphics[width=0.2\linewidth]{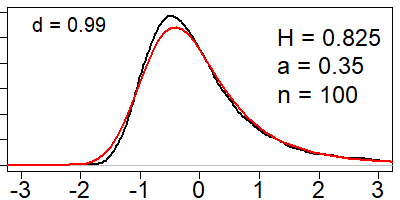}
\includegraphics[width=0.2\linewidth]{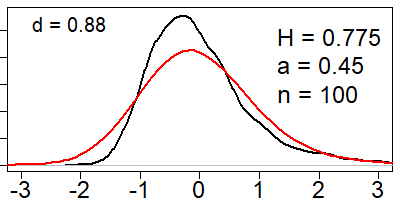}\\
\includegraphics[width=0.2\linewidth]{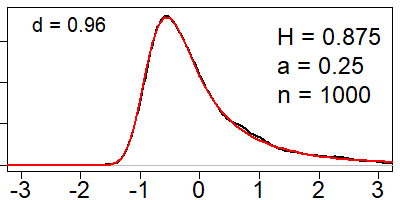}
\includegraphics[width=0.2\linewidth]{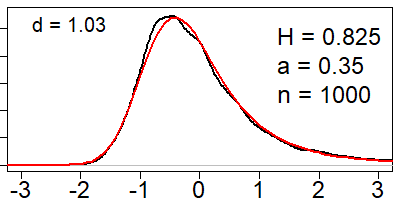}
\includegraphics[width=0.2\linewidth]{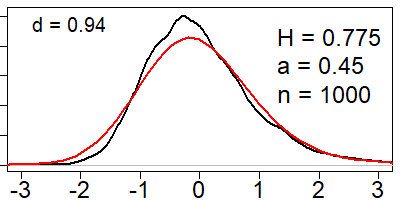}\\
\includegraphics[width=0.2\linewidth]{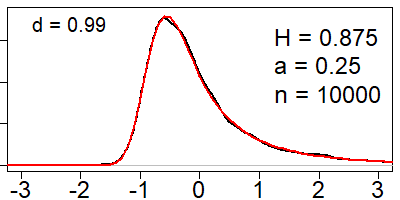}
\includegraphics[width=0.2\linewidth]{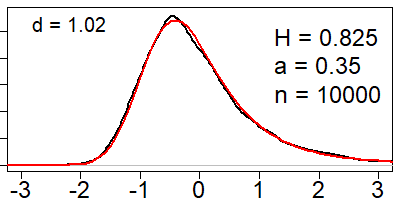}
\includegraphics[width=0.2\linewidth]{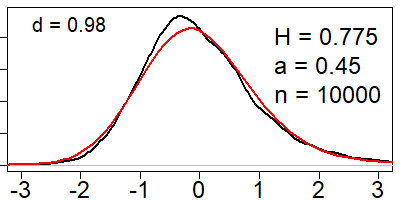}
\end{center}
\caption{The empirical density of $G_{n}$ (black) and the density of the Rosenblatt distribution (red) for various $n$ and $a=0.25$, $0.35$ and $0.45$;
$d$ is the empirical standard deviation of $G_{n}$. }
\label{fig:pat}
\end{figure}

In Figure \ref{fig:pat} we can see that the empirical distribution of $G_{n}$ for small $n$ is rather close to the Rosenblatt distribution for $a\le 0.35$.
We also see that the empirical distribution of $G_{n}$ tends to the Rosenblatt distribution as
$n$ increases for $a=0.45$.

\section{Conclusion}

We studied the Rosenblatt distribution which appears as a limiting distribution of several popular functionals of stationary Gaussian sequences with LRD and, therefore, it is required to construct the confidence intervals.
The Rosenblatt distribution is difficult because it depends on eigenvalues of the Riesz integral operator.
We derived new expressions for the characteristic function of the Rosenblatt distribution, which can be evaluated for any argument.
We obtained the accurate approximation of all eigenvalues that enables easy evaluation of the Rosenblatt distribution. Also, we proposed an efficient algorithm for simulation of stationary Gaussian sequences with the power correlation function $r(t)=1/(1+|t|)^a$ and the Mittag-Leffler correlation function $r(t)=E_{a}(-|t|^a)$.
Finally, we presented Monte-Carlo simulation on how the Rosenblatt distribution appears as a limiting distribution in estimation of several statistics: the mean, the correlation function, sojourn functionals and path roughness.
We note that the Rosenblatt distribution also appears as a limiting distribution of some statistics
for oscillating Gaussian processes 
and in ordinal pattern analysis of stationary Gaussian sequences with LRD.
In addition, the Rosenblatt distribution is the distribution of the Rosenblatt process at unit time, see \cite{pipiras2017long,tudor2023non,ayache2025numerical}.






\ifnum1=1

\section*{Acknowledgments}
 Nikolai Leonenko (NL) would like to thank for support and hospitality during the programme “Fractional Differential Equations” and the programmes “Uncertainly Quantification and Modelling of Material” and "Stochastic systems for anomalous diffusion" in Isaac Newton Institute for Mathematical Sciences, Cambridge.
 These programmes were organized with the support of the Clay Mathematics Institute, of EPSRC (via grants EP/W006227/1 and EP/W00657X/1), of UCL (via the MAPS Visiting Fellowship scheme) and of the Heilbronn Institute for Mathematical Research (for the Sci-Art Contest).
Nikolai Leonenko (NL) was partially supported under the ARC Discovery Grant DP220101680
(Australia), Croatian Scientific Foundation (HRZZ) grant “Scaling in Stochastic Models”
(IP-2022-10-8081), grant FAPESP 22/09201-8 (Brazil) and the Taith Research Mobility
grant (Wales, Cardiff University, 2025).

\fi



\clearpage
\section*{Online supplement. Proof of Theorem 1}

\textbf{(i)}
We note that the logarithm possesses the entire-domain representation
$$
 \ln(x) = 2\sum_{k = 1}^{\infty} \frac{1}{2k-1}\left(\frac{x-1}{x+1}\right)^{2k-1}, ~  x>0,
$$
which is obtained from the Taylor expansion of $\frac{1}{2}\ln(\frac{1+u}{1-u})$ at $u=0$ with substitution $\frac{1+u}{1-u}=x$, see \cite[Eq. 20.19]{spiegel1968mathematical}.
Then we obtain
\begin{eqnarray*}\nonumber
 \ln ( \phi_{LT}(s))&=&- \sum_ {n=1}^\infty \left(
 \frac{\lambda_{a,n}s}{1+\lambda_{a,n}s}  - \lambda_{a,n}s
  + \sum_{k = 2}^{\infty} \frac{1}{2k-1}\left(\frac{\lambda_{a,n}s}{1+\lambda_{a,n}s}\right)^{2k-1}
 \right)\\
 &=& \sum_ {n=1}^\infty \left(
 \frac{(\lambda_{a,n}s)^2}{1+\lambda_{a,n}s}
  - \sum_{k = 2}^{\infty} \frac{1}{2k-1}\left(\frac{\lambda_{a,n}s}{1+\lambda_{a,n}s}\right)^{2k-1}
 \right),
\end{eqnarray*}
which is convenient for numerical computation for all $s>-\frac{1}{2\lambda_{a,1}}$ and
serves as an analytic continuation of the expression \eqref{eq:LTstandard} defined in a neighbourhood of zero.

\textbf{(ii)}
The famous Ramanujan formula is given by
$$
 \frac{1}{\ln(x)}+\frac{1}{1-x}=\sum_{k=1}^\infty \frac{1}{2^{k}(1+x^{2^{-k}})},~  x\in(0,1)\cup(1,\infty),
$$
see \cite{bradley2018concerning}.
We note that $x^{2^{-k}}=x^{(2^{-k})}$ and
$$
 \ln(x)+1-x=(1-x)\ln(x)\left( \frac{1}{\ln(x)}+\frac{1}{1-x} \right)
 =(1-x)\ln(x)\sum_{k=1}^\infty \frac{1}{2^{k}(1+x^{2^{-k}})}.
$$
Then we obtain
\begin{eqnarray*}
 \ln ( \phi_{LT}(s))= \sum_ {n=1}^\infty \lambda_{a,n}s  \ln(1+2\lambda_{a,n}s)
 \sum_{k=1}^\infty \frac{1}{2^{k}(1+(1+2\lambda_{a,n}s)^{2^{-k}})},
\end{eqnarray*}
which is valid for all $s>-\frac{1}{2\lambda_{a,1}}$.

\textbf{(iii)}
The Ramanujan-Bradley expansion is given by
$$
\ln(x)=x-1-\sum_{k=1}^\infty 2^{k-1}(x^{2^{-k}}-1)^2, ~  x>0,
$$
which is similar to the famous Ramanujan formula, see \cite{bradley2017infinite}.
Then we obtain
\begin{eqnarray*}
 \ln ( \phi_{LT}(s))=\frac{1}{2} \sum_ {n=1}^\infty \sum_{k=1}^\infty 2^{k-1}
 \left((1+2\lambda_{a,n}s)^{2^{-k}}-1\right)^2,
\end{eqnarray*}
which is valid for all $s>-\frac{1}{2\lambda_{a,1}}$.

\textbf{(iv)}
The natural logarithm has the integral representation
$$
 \ln(x)=x-1-\int_{1}^x \int_1^t \frac{1}{u^2} du dt.
$$
We can see that
$$
 \ln(x)-x+1=-\int_{1}^x (1-\frac{1}{t})dt=-\int_{0}^{x-1} (1-\frac{1}{u+1})du=-\int_{0}^{x-1} \frac{u}{u+1}du.
$$
For $x=1+2 \lambda_{a,n} s$ we obtain
$$
 \ln(1+2\lambda_{a,n}s )-2\lambda_{a,n}s=
 -\int_{0}^{2 \lambda_{a,n} s} \frac{u}{u+1}du=
  -\int_{0}^{2  s} \frac{\lambda_{a,n}^2u}{\lambda_{a,n} u+1}du.
$$
Then we obtain
$$
 \ln ( \phi_{LT}(s))=\frac{1}{2}\sum_ {n=1}^\infty  \int_{1}^{1+2\lambda_{a,n}s} \int_1^t \frac{1}{u^2} du dt=
 \frac{1}{2}  \int_{0}^{2  s} \sum_ {n=1}^\infty \frac{\lambda_{a,n}^2u}{\lambda_{a,n} u+1}du
 =
   \int_{0}^{ s} \sum_ {n=1}^\infty \frac{2\lambda_{a,n}^2u}{2\lambda_{a,n} u+1}du.
$$

\subsection*{A2. R code for the density of the Rosenblatt distribution}

The following R code contains the  algorithm of computation of
the density of the Rosenblatt distribution with shape parameter $a$.

\begin{lstlisting}


getEig1 = function(a) {
  (1+0.1409*a)*(gamma(1-a)*pi^a)^0.5*(0.5-a)^0.5
}
getEig=function(a, n) {
  S = sqrt((1-2*a)*(1-a)/2)
  Coef = 2*S*gamma(1-a)*sin(pi*a/2)/pi^(1-a)
  SecondCoef = 1.05*a^(5/4)*sqrt(gamma(1-(0.5-a))-1)
  lambda = Coef*n^(a-1) + SecondCoef*n^(a-2.2)
  if(length(n)==1) if (n==1)
    lambda = getEig1(a)
  if(length(a)==1&length(n)>1) if(n[1]==1)
    lambda[1] = getEig1(a)
  return(lambda)
}
IFFT_Fun = function(z, x, a, lam, sigz2) {
  chf = exp(-z^2*sigz2/2)
  for(j in (1:length(lam)))
    chf = chf*exp(-0.5*log(1-2*1i*lam[j]*z)-1i*lam[j]*z)
  return(Re(chf*exp(-1i*z*x)))
}
Ros_Density_Fun = function(x, a, lam, sigz) {
  integrate(IFFT_Fun, -20, 20, x, a, lam, sigz,
            stop.on.error = FALSE)$value/(2*pi)
}
getRosenblattDensity = function(x, a, K) {
  lam = getEig(a,(1:K))
  sigz2 = max(0,1-2*sum(lam^2))
  dens = sapply(x, Ros_Density_Fun, a, lam, sigz2)
  return(dens)
}
a = 0.25
K = 100
x = seq(-3,3,by=0.1)
dens = getRosenblattDensity(x, a, K)
plot(x,dens,type='l')
\end{lstlisting}

\subsection*{A3. R code for simulation of a stationary Gaussian sequence with LRD}

The following R code contains the  algorithm of simulation of
a stationary Gaussian sequence $E_1,E_2,\ldots,E_N$
with zero mean, unit variance and the correlation functions $r(t)=(1+|t|)^{-a}$ and $r(t)=E_a(-|t|^a)$.

\begin{lstlisting}
library(pRSR)
getGausSeqLRD4 = function(N, M, b, lam) {
  X = rep(0,N)
  for(j in (1:M))
  {
   A = exp(-lam[j])
   SEps = sqrt(1-A^2)
   XT = SimulateAR1(N, A)
   X = X + sqrt(b[j])*XT
  }
  return(X)
}
getGausSeqLRD = function(N, a, MLcor=FALSE) {
  gamm = exp((2-a)*a)
  m = ceiling(2/a)+8
  q = if(!MLcor) c(0.98,0.9/gamm^(0:m)) else
              c(0.98,0.9,0.7,0.5/gamm^(0:m))
  Z = (sin(q*a*pi)/sin((1-q)*a*pi))^(1/a)
  Z[Z>10^25] = 10^25
  tau = if(MLcor) Z else qgamma(q, a, 1)
  M = length(q)
  lam = c(tau[1],sqrt(tau[1:(M-1)]*tau[2:M]))
  fun_p = function(x) {
    if(!MLcor) x^(a-1)*exp(-x)/gamma(a) else
      sin(pi*a)/pi*x^(a-1)/(1+2*cos(pi*a)*x^a+x^(2*a)) }
  tau1 = c(4*tau[1], tau)
  b = rep(0, M)
  for(j in (1:M))
    b[j] = integrate(fun_p, tau[j], tau1[j])$value
  X = getGausSeqLRD4(N, M, b, lam)
  return(X)
}
a = 0.35
N = 20000
X = getGausSeqLRD(N, a) # 1/(1+|t|)^a correlation
plot(X, type='l')
acf(X, 100)
X2 = getGausSeqLRD(N, a, TRUE) # E_a(-|t|^a) correlation
\end{lstlisting}



\end{document}